\documentclass[11pt, reqno]{amsart}
\usepackage{amsmath,amsthm,amssymb}

\theoremstyle{plain}
\newtheorem{thm}{Theorem}[section]
\newtheorem{prop}[thm]{Proposition}

\newtheorem{cor}[thm]{Corollary}

\theoremstyle{definition}

\newtheorem{rem}{Remark}
\newtheorem{defn}{Definition}[section]
\newtheorem{eg}[thm]{Example}

\numberwithin{equation}{section}

\newcommand{\bthm}{\begin{thm}}
\newcommand{\ethm}{\end{thm}}
\newcommand{\bprop}{\begin{prop}}
\newcommand{\eprop}{\end{prop}}
\newcommand{\bcor}{\begin{cor}}
\newcommand{\ecor}{\end{cor}}
\newcommand{\bca}{\begin{cases}}
\newcommand{\eca}{\end{cases}}
\newcommand{\brem}{\begin{rem}}
\newcommand{\erem}{\end{rem}}
\newcommand{\bpm}{\begin{pmatrix}}
\newcommand{\epm}{\end{pmatrix}}
\newcommand{\bdefn}{\begin{defn}}
\newcommand{\edefn}{\end{defn}}
\newcommand{\bsub}{\begin{subtitle}}
\newcommand{\esub}{\end{subtitle}}
\newcommand{\ben}{\begin{enumerate}}
\newcommand{\een}{\end{enumerate}}
\newcommand{\beg}{\begin{eg}}
\newcommand{\eeg}{\end{eg}}

\def\ms{\medskip}
\def\bs{\bigskip}

\def\I{{\rm I\/}}
\def\Re{{\rm Re\/}}
\def\Im{{\rm Im\/}}
\def\diag{{\rm diag\/}}

\def\sgn{{\rm sgn\/}}
\def\ti{\tilde}

\def \a {\alpha}

\def \l {\lambda}
\def \L {\Lambda}
\def \n {\,\vert\,}

\def \o {\theta}

\def\R{\mathbb{R} }
\def\C{\mathbb{C}}

\def\Z{\mathbb{Z}}

\def\ce{{\mathcal {E}}}

\def\cm{{\mathcal {M}}}

\def\co{{\mathcal {O}}}

\def \n {\ \vert\ }
\def\tr{{\rm tr}}

\begin{document}


\title[Periodic and homoclinic orbits of the $2+1$ Chiral model]
{Periodic and homoclinic solutions of\\
 the modified $2+1$ Chiral model}
\author{Bo Dai$^*$}\thanks{$^*$Research supported in part by the
AMS Fan Fund\/}
\address{Partner Group of MPI at AMSS\\ Institute of Mathematics\\
Chinese Academy of Sciences\\ Beijing 100080\\ P.~R.~China}
\email{daibo@mail.amss.ac.cn}
\author{ Chuu-Lian Terng$^\dag$}\thanks{$^\dag$Research supported
in  part by NSF Grant DMS-0306446}
\address{Department of Mathematics\\
Northeastern University and University of California at Irvine}
\email{terng@neu.edu}

\ms \hskip 3in \today

\begin{abstract} We use algebraic B\"{a}cklund transformations (BTs) to
construct explicit  solutions of the modified $2+1$ chiral model
from $T^2\times \R$ to $SU(n)$, where $T^2$ is a $2$-torus.
Algebraic BTs are parameterized by $z\in \C$ (poles) and
holomorphic maps $\pi$ from $T^2$ to Gr$(k,\C^n)$.  We apply
B\"acklund transformations with carefully chosen poles and $\pi$'s
to construct infinitely many solutions of  the $2+1$ chiral model
that are (i) doubly periodic in space variables and periodic in
time, i.e., triply periodic, (ii) homoclinic in the sense that the
solution $u$ has the same stationary limit $u_0$ as $t\to
\pm\infty$ and is tangent to a stable linear mode of $u_0$ as
$t\to\infty$ and is tangent to an unstable mode of $u_0$ as $t\to
-\infty$.
\end{abstract}

\maketitle

\section{The $2+1$ Chiral model}

A {\it wave map\/} $J:\R^{2,1}\to SU(n)$ is a critical point of
the functional
$$\ce(J)=\int_{\R^3} ||J^{-1}J_x||^2 + ||J^{-1}J_y||^2 - ||J^{-1}J_t||^2
\ dx dy dt,$$ where $||\xi||^2=-\tr(\xi^2)$, and $x,y,t$ are the
standard space-time variables. The Euler-Lagrange equation of
$\ce$ is
\begin{equation}\label{aa}
(J^{-1}J_t)_t- (J^{-1}J_x)_x - (J^{-1}J_y)_y =0.
\end{equation}
This equation is also called {\it the $2+1$  chiral model\/}.

The {\it Ward equation\/} (or the {\it modified $2+1$ chiral
model\/}) is the following equation for $J:\R^{2,1}\to SU(n)$:
\begin{equation}\label{ward}
(J^{-1}J_t)_t -(J^{-1}J_x)_x -(J^{-1}J_y)_y -[J^{-1}J_t,
J^{-1}J_y]=0.
\end{equation}
This equation is obtained by a dimension reduction and a gauge
fixing of the self-dual Yang-Mills equation on $\R^{2,2}$ (cf.
\cite{W1988}).  We call a solution of the Ward equation a {\it
Ward map\/}.   The Ward
equation is completely integrable and many techniques from
integrable systems can be used to construct explicit solutions.

We consider Ward maps satisfying the doubly periodic boundary
condition in the space variables, i.e., Ward maps from $T^2\times
\R$ to $SU(n)$, where $T^2=S^1\times S^1$.  Using the standard
trick of writing a second order differential equation as a first
order system on the tangent bundle of the phase space, we can view
the Ward equation as a dynamical system on the tangent bundle
$T(C^\infty(T^2, SU(n)))$. The goal of this paper is to construct
periodic and homoclinic orbits of this dynamical system.

A Ward map $J:T^2 \times \R\to SU(n)$ independent of $t$ is a
harmonic map from $T^2$ to $SU(n)$.    Equation for harmonic maps
from $T^2$ to $SU(n)$ is  integrable.  Techniques from integrable
systems were used to construct harmonic maps from $T^2$ to $SU(2)$
by  Hitchin (\cite{Hit90}), and from $T^2$ to $SU(n)$ by    Burstall, Ferus,
Pedit and Pinkall (\cite{BFPP93}).

A Ward map from $S^1\times S^1\times \R$ to $SU(n)$ independent of
the second variable is a wave map from $S^1\times \R$ to $SU(n)$.
Such wave maps were studied by Terng and Uhlenbeck in
\cite{TU2003}.

A solution $u$ of an evolution PDE is {\it homoclinic\/} if  $u$
tends to the same stationary solution $u_0$ as $t\to \pm \infty$
and is tangent to a stable  linear mode of $u_0$ as $t\to +\infty$
and is tangent to an unstable linear mode of $u_0$ as
 $t\to -\infty$. The existence of homoclinic orbits for a finite
dimensional dynamical system indicates the chaotic behavior of the
system (cf. \cite{Hol88}).   It is known that soliton equations in one space and one
time variables (for example, sine-Gordon, KdV, and NLS), viewed as
dynamical systems on certain function spaces,  admit homoclinic
orbits.   Shatah and Strauss \cite{SS1996} proved that there are
homoclinic wave maps from $S^1\times \R$ to $S^2$, and Terng and
Uhlenbeck \cite{TU2003} proved the same result for wave maps from
$S^1\times \R$  to any compact symmetric space. There have been
many works concerning whether homoclinic orbits persist under
small perturbation of these soliton equations in $1$ space and $1$
time variables (cf. \cite{LMSW96, ShaZen03} and references
therein).

One result of this paper is the existence of infinitely many Ward
maps from $T^2\times \R$ to $SU(n)$ that are periodic in time.  In
other words, we prove that there are infinitely many triply
periodic solutions of the Ward equation.  Another result of this
paper is to show that the Ward equation has infinitely many
homoclinic orbits.    We give an outline of our method next.

The $1$-soliton Ward maps from $\R^{2,1}$ to $SU(n)$ can be
constructed as follows (cf. \cite{W1988}). Let $z\in \C\setminus
\R$ be a constant, $V=(v_{ij})$ a meromorphic map from $\C$ to the
space $\cm_{n\times k}^0$ of rank $k$ complex $n\times k$
matrices, $\pi(x,y,t)$ the Hermitian projection of $\C^n$ onto the
subspace spanned by the $k$ columns of $V(w)$, where
 $$w=x+ \frac{(z-z^{-1})y}{2} +\frac{(z+z^{-1})t}{2}.$$
 Let $\pi^\perp= \I-\pi$.    Then
$$\hat J_{z,V}(x,y,t) = \pi^\perp(x,y,t)+\frac{\bar z}{z}
\pi(x,y,t)$$ is a solution of the Ward equation. It has constant
determinant, so we can normalized it to get a Ward map from
$\R^{2,1}$ to $SU(n)$:
$$J_{z,V}(x,y,t) = \left( \frac{z}{\bar z}\right)^{k/n}
\left( \pi^\perp(x,y,t)+\frac{\bar z}{z} \pi(x,y,t)\right).$$
$J_{z,V}$ (or $\hat J_{z,V}$) will be called a {\it Ward
$1$-soliton}. If all entries of $V(w)$ are rational functions in
$w$, then $J_{z,V}$ is a smooth Ward map and is asymptotically
constant as $|(x,y)|\to \infty$.  If all entries of $V$ are
elliptic functions of same periods, then $J_{z,V}$ is a smooth
Ward map from $T^2\times \R$ to $SU(n)$.

Algebraic B\"{a}cklund transformations (BTs) for the Ward equation
were constructed in \cite{DaiTer04}. These are transformations
that generate  new Ward maps from a given Ward map and
$1$-solitons $J_{z, V}$ by a simple algebraic method.

 We apply algebraic
B\"acklund transformations repeatedly to  $1$-solitons associated
to elliptic functions to construct infinitely many triply periodic
Ward maps to $SU(n)$

Note that if  the image of $J$ lies  in an abelian subgroup of
$SU(n)$, then the Ward equation for $J$ becomes the linear wave
equation.  For example, let $m$ be an integer, and $a=\diag(im,
-im)$.  Then
 $$J_0(x,y,t)= \exp(-(x+y)a)$$
 is a doubly periodic, stationary Ward map, whose image lies in $SO(2)$.

 We apply algebraic BTs $2k$ times  to   $J_0$ with carefully chosen
poles and projections to construct homoclinic Ward maps $J_{2k}$
from $T^2\times \R$ to $SU(n)$, and prove that $J_{2k}$ tends to
$(-1)^kJ_0$ as $|t|\to\infty$ and $J_{2k}$ are homoclinic.

This paper is organized as follows:  We review the Lax pair  and
algebraic B\"acklund transformations for the Ward equation in
section 2, and  use elliptic functions to construct triply
periodic Ward maps  in section 3.  In the last section,  we
construct (i) homoclinic Ward maps from $T^2\times \R$ to $SU(n)$
that tend to stationary solutions, (ii)  homoclinic Ward maps that
tend to periodic solutions.

\bs
\section{Extended Ward maps and B\"acklund transformations}

 The Ward equation is integrable in the sense that it can be
written as the compatibility condition for a system of linear
equations involving a spectral parameter $\lambda\in\C$. In fact,
we have the following theorem (cf. \cite{W1988})

 \bthm  \label{dt}
Let $J:\R^{2,1}\to SU(n)$ be a Ward map, $dx^2+dy^2-dt^2$ be the
Lorentzian metric on $\R^{2,1}$,
\begin{equation} \label{ai}
u=\frac{t+y}{2}, \quad v=\frac{t-y}{2},
\end{equation}
 $A=J^{-1}J_u$, and $B=J^{-1}J_x$.  Then the following linear PDE
system is solvable for $\psi:\R^{2,1}\times\C \to GL(n,\C)$:
\begin{equation} \label{eq:auxiliary}
\bca
(\lambda\partial_x -\partial_u)\psi =A\psi, \\
(\lambda\partial_v -\partial_x)\psi =B\psi. \eca
\end{equation}

Conversely, suppose $\co$ is an open subset of $0$ in $\C$ and $
\psi:\R^{2,1}\times \co\to GL(n,\C)$ is a smooth map so that
$$A:=(\l\psi_x-\psi_u)\psi^{-1}, \quad B:= (\l\psi_v-\psi_x)\psi^{-1}$$
are independent of $\l\in \co$ and $\psi$ satisfies the
$U(n)$-reality condition
\begin{equation}
\psi (x,u,v,\bar{\lambda})^* \psi (x,u,v,\lambda)={\rm I\/},
\label{eq:reality}
\end{equation}
Then
$$J(x,y,t)= \psi(x,y,t,0)^{-1}$$
is a smooth solution of the Ward equation and $J^{-1}J_u=A$ and
$J^{-1}J_x=B$. \ethm

A solution $\psi(x,y,t,\l)$ of \eqref{eq:auxiliary} that satisfies
the $U(n)$-reality condition \eqref{eq:reality} is called an {\it
extended Ward map\/} and $J=\psi(\cdots, 0)^{-1}$ the associated
Ward map.

Given $z\in \C$ and a Hermitian projection $\pi$ of $\C^n$, let
$$h_{z,\pi}(\l)= \pi^\perp+ \frac{\l-z}{\l-\bar z}\pi = \I +
\frac{\bar z-z}{\l-\bar z}\pi,$$ where $\pi^\perp= \I-\pi$.  A
direct computation implies that $h_{z,\pi}$ satisfies the
$U(n)$-reality condition \eqref{eq:reality}.

Let $V=(v_{ij}):\C\to \cm_{n\times k}^0(\C)$ be a meromorphic map,
and $\pi(x,y,t)$ the Hermitian projection onto the subspace
spanned by the columns of $V(w)$, where
$$w=x+zu + z^{-1}v,$$
and $u,v$ are the light cone coordinates in the $yt$-plane defined by \eqref{ai}. Since the entries of $V$ are
meromorphic functions, the projection $\pi$ is smooth on
$\R^{2,1}$.  Set
$$\psi(x,y,t,\l)= h_{z, \pi(x,y,t)}(\l)= \pi^\perp (x,y,t)
+\frac{\l-z}{\l-\bar z} \, \pi (x,y,t)$$ A direct computation
implies that both $(\l\psi_x-\psi_u)\psi^{-1}$ and
$(\l\psi_v-\psi_x)\psi^{-1}$ are independent of $\l$.  By Theorem
\ref{dt}, $\psi$ is an extended solution of the Ward equation and
the associated Ward map is the $1$-soliton
$$J_{z,V}(x,y,t)=\left( \frac{z}{\bar z}\right)^{k/n}
\psi(x,y,t, 0)^{-1} =\left( \frac{z}{\bar z}\right)^{k/n} \left(
\pi^\perp(x,y,t) +\frac{\bar z}{z} \pi(x,y,t)\right),$$ where
$\left( \frac{z}{\bar z}\right)^{k/n}$ is a normalizing constant
to make $\det (J_{z,V})=1$.

The $1$-soliton Ward map $J_{z,V}$ is a travelling wave because
$$w=x+zu+z^{-1}v =(x-v_1 t) +k_1(y-v_2 t) +ik_2(y-v_2 t),$$
where $z=re^{i\o}$, $v_1=-\frac{2r \cos\theta}{1+r^2}$,
$v_2=\frac{1-r^2}{1+r^2}$, and $k_1 +ik_2=(z -z^{-1})/2$. Thus
$J_{z,V}$ is a travelling wave with constant velocity $\vec{v}
=(-\frac{2r \cos\theta}{1+r^2}, \frac{1-r^2}{1+r^2})$ on the
$xy$-plane. In particular, $J_{i,V}$ is a stationary Ward map,
i.e., a harmonic map from $\C$ to $SU(n)$.

The following theorem was proved in \cite{DaiTer04}, which gives
an algebraic method to produce new extended Ward maps from a given
one.

\begin{thm}[B\"{a}cklund transformation]
\label{thm:BT} Let $\psi(x,y,t,\l)$ be an extended solution of the
Ward equation  and $J= \psi(\cdots, 0)^{-1}$ the associated Ward
map from $\R^{2,1}$ to $SU(n)$. Choose $z\in {\mathbb C}\setminus
{\mathbb R}$ such that $\psi(x,y,t,\lambda)$ is holomorphic and
non-degenerate at $\lambda =z$. Let $h_{z,\pi(x,y,t)}(\lambda)$ be
an extended $1$-soliton solution, and $\ti\pi(x,y,t)$ the
Hermitian projection of $\C^n$ onto
$$\psi(x,y,t,z)\Im(\pi(x,y,t)).$$ Then
$$\psi_1 (x,y,t,\lambda) = h_{z,\tilde{\pi}(x,y,t)}(\lambda)
\psi(x,,y,t,\lambda)$$ is a new extended solution to the linear
system (\ref{eq:auxiliary}) with $$(A,B)\to (A +(\bar{z}-z)
\tilde{\pi}_x, B + (\bar{z}-z) \tilde{\pi}_v),$$ and the new Ward
map is
$$J_1(x,y,t) = \left( \frac{z}{\bar z}\right)^{k/n}
J(x,y,t)\, \left( \frac{\bar{z}}{z} \tilde{\pi} (x,y,t)
+{\tilde{\pi}}^\perp (x,y,t) \right).$$
\end{thm}

We will  denote $\psi_1=h_{z,\pi} * \psi$ and $J_1=h_{z,\pi} *J$,
the B\"{a}cklund transformation generated by $h_{z,\pi}$.

\bs

\section{Periodic Ward maps from $T^2\times \R$ to $SU(n)$}

We use algebraic BTs to construct Ward maps into $SU(n)$ that are either
doubly periodic in space variables or triply periodic.

First we construct $1$-soliton Ward maps that are doubly periodic.
Let $z=re^{i\o}$,
$$w(x,y,t)= x+zu+z^{-1}v= x+\frac{z-z^{-1}}{2} \ y +
\frac{z+z^{-1}}{2} \ t,$$ and $\a= a+ib$.
  A direct computation shows that
$$w(x+a-\frac{k_1}{k_2}\ b,\  y+ \frac{b}{k_2},\  t)= w(x,y,t)+\a,$$
where $k_1+ik_2= \frac{z-z^{-1}}{2}$. Let $f:\C\to \C$ be a
meromorphic function such that $f(w+\a)=f(w)$ (i.e., periodic with
period $\a$),  and
$$g(x,y,t)= f(w)= f(x+zu+z^{-1}v).$$
A direct computation shows that
$$g(x,y,t)= g(x+a-\frac{k_1}{k_2}b, y+\frac{b}{k_2} ,t).$$
 Hence
 \ben
 \item if $\alpha =2\pi$,
then $g$ is $2\pi$-periodic in $x$;
 \item
 if $\alpha =2\pi (k_1
+ik_2) =\pi(z-z^{-1})$, then $g$ is $2\pi$-periodic in $y$.
 \een

This shows that if each entry $v_{ij}$ of the meromorphic map
$V=(v_{ij}):\C\to \cm_{n\times k}^0(\C)$ satisfies $v_{ij}(w
+2\pi) =v_{ij}(w +\pi(z-z^{-1})) =v_{ij}(w)$, i.e. an elliptic
function with periods $2\pi$ and $\pi(z-z^{-1})$,  then the
$1$-soliton $J_{z,V}$ is a doubly periodic Ward map with respect
to the lattice $2\pi(\Z\times \Z)$. An example of elliptic
function is the well-known Weierstrass $\wp$-function
$$\wp (w) =\frac{1}{w^2} + \sum_{\gamma\in\Lambda
\setminus\{ 0\}} \left( \frac{1}{(w-\gamma)^2}
-\frac{1}{\gamma^2}\right),$$ where $\Lambda$ is the lattice in
${\mathbb C}$ generated by $2\pi$ and $\pi(z-z^{-1})$. Other
elliptic functions can be generated by Weierstrass $\wp$-functions
and Jacobi elliptic functions. It is clear that $J_{z,V}$ is time
periodic if and only if the ratio of the velocity, $v_1/v_2
=(-2r\cos\theta)/(1-r^2)$, is rational.

Similar computation implies that given any rank $2$ lattice $\L$
of $\C$ there are $1$-soliton Ward maps from $\C/\L\times \R$ to
$SU(n)$.  Moreover, some of these $1$-solitons are periodic in
time, i.e., triply periodic.  In particular, we get

\bthm \label{dp} Let $\tau=c_1+ ic_2$ with $c_2\not=0$,  $\L= \Z
2\pi + \Z \tau$, $z=re^{i\o}$ a constant, and $a+ib= c_1+
\frac{z-z^{-1}}{2} c_2$.   If each entry $v_{ij}$ of the
meromorphic map $V:\C\to \cm_{n\times k}^0(\C)$ is an elliptic
function with periods $2\pi$ and  $a+ib$, then the extended
$1$-soliton solution $h_{z, \pi}$ is doubly periodic with periods
$2\pi$ and $\tau$ and the associated $1$-soliton
$$J_{z, V}=e^{i2k\o /n} h_{z,\pi}(0)^{-1}=
e^{i2k\o /n} (\pi^\perp+ e^{-2i\o}\pi)$$ is a Ward map from
$\C/\L\times \R$ to $SU(n)$, where $e^{i2k\o /n}$ is a normalizing
constant, and $\pi(x,y,t)$ is the projection onto the subspace
spanned by the columns of $V(x+zu+ z^{-1}v)$. Moreover, \ben \item
if $r\not=1$ and there exist integers $m_1, m_2$ such that
$$\frac{2\cos \o}{r-r^{-1}}= \frac{2\pi m_1 + m_2 c_1}{m_2 c_2},$$ then
$J_{z, V}$ is periodic in time with period
$T=\frac{m_2c_2(r+r^{-1})}{r-r^{-1}}$, \item if $r=1$ and
$\cos\o\not=0$, then $J_{z,V}$ is periodic in time with period $T=
\frac{2\pi}{\cos\o}$. \een \ethm

In the rest of the section we consider only the square torus.  We
will construct $k$-soliton Ward maps from $T^2\times \R$ to
$SU(n)$ that are also time periodic.  To do this,  we define
$${\mathcal Z}= \{ z=re^{i\theta} \in {\mathbb C}\setminus
{\mathbb R}\, \  \vert\  z=e^{i\o}\not=\pm i \ {\rm or\/}\
\cos\theta /(r-r^{-1})\in {\mathbb Q} \},$$ where ${\mathbb Q}$
denotes the set of rational numbers. We have seen that for each
$z\in {\mathcal Z}$, we can construct time periodic $1$-solitons
to the Ward equation. Moreover, the time period $T$ depends on $z$
only.  In fact, the period function $T: {\mathcal Z}\to \R$ is
defined as follows:
\begin{equation}\label{dq}
T(z)=\bca \frac{2\pi}{\cos \o}, & {\rm if \ } z=e^{i\o}\not=\pm i,\\
\frac{2\pi m_2(r+r^{-1})}{r-r^{-1}}, & {\rm if\ } z=re^{i\o}, \
\frac{2\cos\o}{r-r^{-1}}= \frac{m_1}{m_2}. \eca
\end{equation}
Apply B\"{a}cklund transformations repeatedly with some rational
conditions on the poles $z_1, \ldots, z_m$ to get the following:

\begin{thm} Let $\{ z_1,\dots,z_m \}$ be a set of finite points in
${\mathcal Z}$ such that $z_i \not= z_j,\bar{z}_j$ for all $i\not=
j$, and $h_{z_i,\pi_i}(\lambda)$ extended $1$-soliton solutions
leading to doubly periodic Ward maps, where $i,j=1,\cdots m$. Let
$T_i= T(z_i)$ be the time period  defined in \eqref{dq}. Let $J_1$
be the Ward map associated to $h_{z_1,\pi_1}$, i.e., $J_1= h_{z_1,
\pi_1}(0)^{-1}$. Let $J_m$ be the Ward map obtained by applying
$m-1$ B\"acklund transformations to $J_1$,
\begin{equation} J_m =h_{z_m,\pi_m}*(\cdots *(h_{z_2,\pi_2}*J_1)
\cdots ).
\end{equation}
 If $T_j/T_1$ are rational numbers for all $2\leq j\leq m$, then
$J_m$ is a Ward map from $T^2\times \R$ to $SU(n)$ and is periodic
in time.  In other words, $J_m$ is a triply periodic solution of
the Ward equation.
\end{thm}

\begin{proof}
We prove the two-soliton case. By Theorem~\ref{thm:BT}, we have
$$h_{z_2,\pi_2}*h_{z_1,\pi_1} =h_{z_2,\tilde{\pi}_2}
h_{z_1,\pi_1},$$ where $\Im\tilde{\pi}_2 =h_{z_1,\pi_1}(z_2) \Im
\pi_2 =({\rm I\/} +\frac{\bar z_1-z_1}{z_2- \bar z_1} \pi_1) \Im
\pi_2$.  Note that $\tilde{\pi}_2$ is periodic in time because
$\pi_1$ and $\pi_2$ are time periodic and $T_2/T_1$ is rational.
Thus we see that $h_{z_2,\pi_2}*h_{z_1,\pi_1}$ is time periodic,
and so is the associated Ward map.  The general case can be proved by
induction.
\end{proof}

\bs
\section{Homoclinic Ward maps}

It is known that solutions of the sine-Gordon equation (SGE) 
$$q_{tt}-q_{xx}= \sin q$$
give rise to wave maps from $\R^{1,1}$ to $S^2$.   Breather solutions
are $2$-soliton solutions of the SGE that are periodic in the $x$ variable.
Shatah and Strauss proved in \cite{SS1996} that wave maps from
$S^1\times \R$ to $S^2$ corresponding to breather solutions of the
sine-Gordon equation are homoclinic wave maps. Applying
B\"{a}cklund transformation $2k$-times with carefully placed
poles, Terng and Uhlenbeck constructed $2k$-soliton solutions for
the sine-Gordon equation that are periodic in the space variable,
and showed that the corresponding wave maps from $S^1\times \R$ to
$S^2$ are also homoclinic. More generally they proved that there
are homoclinic wave maps from $S^1\times \R$ into any compact
symmetric space \cite{TU2003}.

In this section, we apply B\"{a}cklund transformations  with
carefully chosen poles and Hermitian projections even times to
certain stationary wave map into $SO(2)$ to construct homoclinic
Ward maps from $T^2\times \R$ to $SU(n)$. To make the construction
more illuminating, we will work on the $SU(2)$ model. The $SU(n)$
model is similar.

Let $m>0$ be an integer, and $a=\text{diag} (im, -im)\in su(2)$.
It is easy to check that
\begin{equation} \psi(\l)(x,y,t)=\psi(x,y,t,\l)
=e^{((1-\lambda)x +(1+\lambda-\lambda^2)u
-v)a}. \label{trivialpsi} \end{equation} is an extended Ward map.
So
$$J_0(x,y,t)=\psi (x,y,t,0)^{-1} =e^{-(x+u-v)a} =e^{-(x+y)a}$$ is a
stationary Ward map, which is doubly periodic in the space
variables. Note that $J_0$  is a harmonic map from $T^2$ to $SO(2)$.

Next we compute the linearization of the Ward equation at the
stationary solution $J_0=e^{-(x+y)a}$, as well as its stable and
unstable subspaces.  Let ${\mathcal M}=C^\infty (T^2 \times
{\mathbb R},SU(2))$. Then we can give a natural trivialization of
the tangent bundle $T\cm$ as follows. Given a curve $\gamma:
(-\epsilon,\epsilon) \to {\mathcal M}$ with $\gamma(0) =J$, we
identify the tangent vector $\gamma'(0)$ as
$$(\gamma(0), \gamma(0)^{-1} \gamma'(0)) =(J, J^{-1} \delta J).$$
This identifies $T{\mathcal M}= {\mathcal M}\times C^\infty (T^2
\times {\mathbb R},su(2))$.

Set $J^{-1}\delta J =\eta$. Compute directly to get
\begin{eqnarray*} \delta (J^{-1}J_x) &=& -(J^{-1}\delta J)
J^{-1}J_x +J^{-1} (\delta J)_x
\\ &=& -\eta (J^{-1}J_x) +J^{-1} (J\eta)_x =-\eta (J^{-1}J_x)
+J^{-1} (J_x\eta +J\eta_x) \\ &=& \eta_x +[J^{-1}J_x,\eta].
\end{eqnarray*} The computation for $\delta (J^{-1}J_y)$ and $\delta
(J^{-1}J_t)$ is similar. So the linearization of the Ward equation
at $J_0=e^{-(x+y)a}$ is:
\begin{eqnarray}
& &(\eta_t +[J^{-1}J_t,\eta])_t - (\eta_x +[J^{-1}J_x,\eta])_x
-(\eta_y +[J^{-1}J_y,\eta])_y \notag \\ & & \quad - [\eta_t
+[J^{-1}J_t,\eta],J^{-1}J_y] -[J^{-1}J_t, \eta_y
+[J^{-1}J_y,\eta]] \notag \\ &=& \eta_{tt}- \eta_{xx} -\eta_{yy}
+[a, \eta_x +\eta_y -\eta_t] =0. \label{linearize}
\end{eqnarray} We note that the linearization at
$J=-e^{-(x+y)a}$ is the same one. Write (\ref{linearize}) in terms
of entries $\eta =\begin{pmatrix} ir & \xi \\ -\bar{\xi} & -ir
\end{pmatrix}$ to get \begin{eqnarray} \left\{
\begin{array}{l}
r_{tt}-r_{xx}-r_{yy} =0, \\
\xi_{tt}-\xi_{xx} -\xi_{yy} +2im (\xi_x +\xi_y -\xi_t)=0.
\end{array} \right. \label{eqs:rxi} \end{eqnarray}
This system is linear with constant coefficients, so it can be
solved by Fourier series. Let $$\xi =\sum_{j,l\in{\mathbb Z}}
b_{jl}(t) e^{i(jx+ly)}$$ be the Fourier series expansion of $\xi$.
Then by (\ref{eqs:rxi}.2), we have $$b_{jl}^{\prime\prime}
-2imb_{jl}^\prime +(j^2+l^2 -2m(j+l)) b_{jl} =0,$$ where $\prime$
means differentiation with respect to $t$. Its auxiliary equation
is $$\gamma^2 -2im\gamma +j^2+l^2 -2m(j+l)=0.$$ It has roots
$$\gamma =im \pm \sqrt{m^2-(j-m)^2 -(l-m)^2}.$$ Stable (unstable
respectively) modes come from $\Re(\gamma)  <0$ ($\Re(\gamma)> 0$
respectively). So for $(j,l) \in {\mathbb Z}^2$ with $(j-m)^2
+(l-m)^2  <m^2$, there are stable and unstable modes corresponding
to roots $im \mp \sqrt{m^2-(j-m)^2 -(l-m)^2}$ respectively.
Similar computation shows that the auxiliary equation for
(\ref{eqs:rxi}.1) has only purely imaginary roots. So the above
computation gives

\begin{prop} \label{ak}
Let $a=\text{diag}(im,-im)$ and $J=\pm e^{-(x+y)a}$, where  $m>0$
is an integer.  Let
$$B\Z_m=\{(j,l)\in \Z^2 \n  (j-m)^2 +(l-m)^2  <m^2\}.$$
 Then:
 \ben
 \item The unstable subspace of the linearization of the
Ward equation at $J$ is
 $$\bigoplus \{W_{jl}^+\n (j,l) \in B\Z_m\},$$
 where  $W_{jl}^+ $ is
spanned by
$$\eta_{jl}^+ (c) =e^{\sqrt{m^2-(j-m)^2
-(l-m)^2}\ t} \begin{pmatrix} 0 & ce^{i(jx+ly+mt)} \\ -\bar{c}
e^{-i(jx+ly+mt)} & 0 \end{pmatrix}$$ with constant $c\in \C$.
 \item
 The stable subspace at $J$ is
 $$\bigoplus \{W_{jl}^-\n (j,l)\in B\Z_m\},$$
where $ W_{jl}^-$ is spanned by
$$\eta_{jl}^-
(c) =e^{-\sqrt{m^2-(j-m)^2 -(l-m)^2}\ t} \begin{pmatrix} 0 &
ce^{i(jx+ly+mt)} \\ -\bar{c} e^{-i(jx+ly+mt)} & 0 \end{pmatrix}$$
with $c\in \C$.
 \een
\end{prop}

Let $z=re^{i\theta}\in {\mathbb C} \setminus {\mathbb R}$,
$f(w)$ a meromorphic function on $\C$, $q(w) =\begin{pmatrix} 1\\
f(w) \end{pmatrix}$,  $w =x +zu +z^{-1}v$, and $\pi(x,y,t)$ the
Hermitian projection of ${\mathbb C}^2$ onto ${\mathbb C} q(w)$.
Let $\psi$ be the extended solution given by (\ref{trivialpsi})
and $J_0=\psi^{-1}|_{\l=0}= e^{-a(x+y)}$ the associated Ward map.
Consider the B\"{a}cklund transformation $h_{z,\pi} *\psi$. We
will find conditions on $z$ and $f(w)$ so that $h_{z,\pi} *J_0$ is
doubly periodic in space variables. By Theorem \ref{thm:BT}, we
have
\begin{equation}\label{aj}
\psi_1 = h_{z,\pi} *\psi =h_{z, \tilde{\pi}} \psi,
\end{equation}
where $\text{Im} \tilde{\pi}  ={\mathbb C} \tilde{q}$ and 
 \begin{eqnarray*} 
 \tilde{q}(x,y,t) &=& \psi(z) q(w) \\
&=& e^{((1-z)x +(1+z-z^2)u-v)a} \begin{pmatrix} 1\\ f(w)
\end{pmatrix} \\ &\sim& \begin{pmatrix} 1\\  \\  e^{2im((z-1)x
+(z^2 -z-1)u+v)} f(w) \end{pmatrix}.
  \end{eqnarray*} 
  Here
``$q_1\sim q_2$'' means $\C q_1=\C q_2$. From the formula
$$J_1
=h_{z,\pi} *J_0 =J_0\frac{1}{\vert z\vert} (\bar{z}\tilde{\pi}
+z\tilde{\pi}^\perp),$$
 we see that it is doubly periodic if and
only if $\tilde{q}$ is. For this purpose, we try the following
form of
$$f(w) =e^{2im (\alpha -z)w},$$ where $\alpha \in
{\mathbb C}$ is a constant. Substitute this into $\ti q(x,y,t)$ to
get
\begin{eqnarray*} e^{2im((z-1)x +(z^2 -z-1)u+v)} f(w) =
e^{2im((\alpha -1)x +((\alpha -1)z-1)u+ \alpha z^{-1}v)} \\
=e^{im(2(\alpha -1)x +((\alpha -1)z -\alpha z^{-1}-1)y +((\alpha
-1)z +\alpha z^{-1}-1)t)}.\end{eqnarray*} It is doubly periodic in
$x$ and $y$ with period $2\pi$ if and only if
\begin{equation}
\begin{cases}
 2m(\alpha -1):= -j \in {\mathbb Z}, & \\
 m((\alpha -1)z -\alpha z^{-1}-1) :=-l \in {\mathbb Z}. &
\end{cases} \label{jl:cond} \end{equation} From
(\ref{jl:cond}.1), we have $\alpha =\frac{2m -j}{2m}$. Compute the
imaginary part of (\ref{jl:cond}.2) to get
$$(\alpha -1)r\sin \theta +\alpha r^{-1} \sin \theta =0.$$ Since
$r>0$, we see $0< \alpha  <1$. This implies that $0<j<2m$, and
$r=\sqrt{\frac{2m-j}{j}}$. By (\ref{jl:cond}.2) again, we have
$$ m((\alpha -1)z -\alpha z^{-1}-1) =-\sqrt{j(2m -j)}\cos \theta
-m=-l.$$ It follows that $$ \sqrt{j(2m -j)}\cos \theta =l-m.$$
Hence $l$ must satisfy
\begin{equation} \vert l-m \vert  < \sqrt{j(2m -j)}, \label{lcond}
\end{equation}
and $\cos \theta =\frac{l-m}{\sqrt{j(2m -j)}}$. It is easy to
verify that the conditions for $(j,l)$ are equivalent to $(j,l)\in
{\mathbb Z}^2$, $(j-m)^2 +(l-m)^2<m^2$, i.e., $(j,l)\in B\Z_m$.
Therefore if we choose the following data: $(j,l)\in B\Z_m$,
$z=re^{i\theta}$ with $r=\sqrt{\frac{2m-j}{j}}$, $\cos \theta
=\frac{l-m}{\sqrt{j(2m -j)}}$, $\sin\theta>0$, $\alpha =\frac{2m
-j}{2m}$, then $\text{Im} \tilde{\pi}(x,y,t) ={\mathbb C}
\tilde{q}(x,y,t)$, where
$$ \tilde{q}(x,y,t) =\begin{pmatrix} 1\\ \\ e^{\sqrt{m^2 -
(j-m)^2 -(l-m)^2}\ t}\ e^{-i(jx+ly+mt)}
\end{pmatrix}$$ is doubly periodic in $x$ and $y$. It follows
that
\begin{equation*} \tilde{\pi}(x,y,t) = \frac{1}{1+ e^{2A}}
\begin{pmatrix} 1 & e^A e^{i(jx+ly+mt)} \\ e^A e^{-i(jx+ly+mt)} & e^{2A}
\end{pmatrix}.
\end{equation*}
where $A= {\sqrt{m^2 - (j-m)^2 -(l-m)^2}\ t}$. Therefore we obtain
the following Ward map from $T^2 \times\R$ to $SU(2)$:
\begin{eqnarray} J_1 &=&
e^{-(x+y)a}\frac{1}{\vert z\vert}
(\bar{z} \tilde{\pi}(x,y,t) +z\tilde{\pi}^\perp(x,y,t)) \notag \\
&=& e^{-(x+y)a} (e^{-i\theta} \tilde{\pi}(x,y,t) +e^{i\theta}
\tilde{\pi}^\perp(x,y,t)). \end{eqnarray} We now analyze the
asymptotic behavior of $J_1$ as $t \to \pm\infty$. It is easy to
see that
$$\tilde{\pi} \to
\begin{pmatrix} 0 & 0\\ 0 & 1 \end{pmatrix} \quad \text{as}\
t\to +\infty,$$ and $$\tilde{\pi} \to
\begin{pmatrix} 1 & 0\\ 0 & 0
\end{pmatrix} \quad \text{as}\ t\to -\infty.$$ So $$J_1 \to e^{-(x+y)a}
\begin{pmatrix} e^{i\theta} & 0\\ 0 & e^{-i\theta}
\end{pmatrix} \quad \text{as}\ t\to +\infty,$$ and $$J_1 \to
e^{-(x+y)a} \begin{pmatrix} e^{-i\theta} & 0\\ 0 & e^{i\theta}
\end{pmatrix}\quad \text{as}\ t\to -\infty.$$ Thus $J_1: T^2
\times {\mathbb R}\to SU(2)$ is a heteroclinic Ward map. To
construct homoclinic maps, we apply B\"{a}cklund transformation
again.

Choose $z_2 =-\bar{z}$, and $\pi_2 (x,y,t)$ Hermitian projection
of ${\mathbb
C}^2$ onto ${\mathbb C}q_2$, where $q_2 =\begin{pmatrix} 1 \\
f_2(w_2) \end{pmatrix}$, $f_2(w_2) =e^{2im (\alpha +\bar{z})w_2}$,
$\a= \frac{2m-j}{2m}$,  and $w_2=x -\bar{z}u -\bar{z}^{-1}v$. Now
apply B\"{a}cklund transformation to $\psi_1$ (defined by
\eqref{aj}) generated by $h_{-\bar{z}, \pi_2(x,y,t)}$ to get
$$\psi_2 =h_{-\bar{z}, \pi_2(x,y,t)} *\psi_1 =h_{-\bar{z},
\tilde{\pi}_2(x,y,t)} \psi_1,$$ where
$\text{Im}\tilde{\pi}_2(x,y,t) ={\mathbb C} \tilde{q}_2 (x,y,t)$
and
\begin{eqnarray} \tilde{q}_2 &=& \psi_1(-\bar{z}) q_2
=h_{z,\tilde{\pi}(x,y,t)}(-\bar{z}) \psi(-\bar{z})
\begin{pmatrix} 1 \\ f_2(w_2) \end{pmatrix} \notag \\ &\sim&
\left({\rm I\/}+\frac{\bar{z}-z}{-2\bar{z}} \tilde{\pi}(x,y,t)
\right) \begin{pmatrix} 1\\ e^A e^{-i(jx+(2m-l)y+mt)}
\end{pmatrix}. \label{as}
\end{eqnarray}
 Here $A= {\sqrt{m^2 - (j-m)^2 -(l-m)^2}\ t}$.  Note that $\ti
q(x,y,t)$ is doubly periodic in $x$ and $y$.  Therefore
\begin{eqnarray*} &J_2  =& J_0\frac{1}{\vert
z\vert} (\bar{z} \tilde{\pi} +z\tilde{\pi}^\perp) \frac{1}{\vert
z_2\vert} (\bar{z}_2 \tilde{\pi}_2 +z_2\tilde{\pi}_2^\perp) \\ &=&
J_0 (e^{-i\theta}\tilde{\pi} +e^{i\theta} \tilde{\pi}^\perp)
(-e^{i\theta} \tilde{\pi}_2 -e^{-i\theta} \tilde{\pi}_2^\perp)
\end{eqnarray*} is a Ward map from $T^2\times\R$ to $SU(2)$,
where $J_0 =e^{-(x+y)a}$.

Next we study the asymptotic behavior of $J_2$. First look at the
behavior of $J_2$ as $t\to -\infty$. Set $$\xi =e^{\sqrt{m^2 -
(j-m)^2 -(l-m)^2}\ t},\ f_1 =e^{-i(jx+ly+mt)},\ f_2
=e^{-i(jx+(2m-l)y+mt)}.
$$ Then $\lim_{t\to -\infty}\xi=0$. Write $$\tilde{q}_1 =\ti{q}=
\begin{pmatrix} 1\\ \xi f_1 \end{pmatrix}  = \begin{pmatrix} 1\\
0 \end{pmatrix} +\xi f_1 \begin{pmatrix} 0 \\ 1 \end{pmatrix}.$$
So the projection $\tilde{\pi}_1=\ti\pi$ onto $\C \tilde{q}_1$ is
$$\tilde{\pi}_1 =\begin{pmatrix} 1 & 0 \\ 0 & 0 \end{pmatrix} +
\xi \begin{pmatrix} 0 & \bar{f}_1 \\ f_1 & 0 \end{pmatrix}
+O(\xi^2).$$  Write $\alpha_1 =\frac{z-\bar{z}}{2\bar{z}}$. Then
by \eqref{as} we have
\begin{eqnarray*} \tilde{q}_2 &=&
({\rm I\/}+\alpha_1 \tilde{\pi}_1)
\begin{pmatrix} 1\\ \xi f_2 \end{pmatrix} \\ &=&
({\rm I\/}+\alpha_1 \tilde{\pi}_1) \left(\begin{pmatrix} 1\\
0 \end{pmatrix} +\xi f_2 \begin{pmatrix} 0 \\ 1 \end{pmatrix} \right)\\
&=& \begin{pmatrix} 1+\alpha_1 \\ 0 \end{pmatrix} + \alpha_1 \xi
f_1
\begin{pmatrix} 0 \\ 1 \end{pmatrix} + \xi f_2 \begin{pmatrix}
0 \\ 1 \end{pmatrix} +O(\xi^2) \\
&\sim& \begin{pmatrix} 1\\ 0 \end{pmatrix} +\beta_1 \xi f_1
\begin{pmatrix} 0 \\ 1 \end{pmatrix} +\beta_2 \xi f_2
\begin{pmatrix} 0 \\ 1 \end{pmatrix} +O(\xi^2), \end{eqnarray*}
where $\beta_1 =\frac{\alpha_1}{1+\alpha_1}$, $\beta_2
=\frac{1}{1+\alpha_1}$. So the projection $\tilde{\pi}_2$ onto
${\mathbb C}\tilde{q}_2$ is
$$\tilde{\pi}_2  = \begin{pmatrix} 1 & 0 \\ 0 & 0 \end{pmatrix}
+\xi \begin{pmatrix} 0 & \bar{\beta}_1\bar{f}_1 \\ \beta_1 f_1 & 0
\end{pmatrix} +\xi \begin{pmatrix} 0 & \bar{\beta}_2\bar{f}_2 \\
\beta_2 f_2 & 0 \end{pmatrix} +O(\xi^2).$$ From the above
computation, we see $$\lim_{t\to -\infty} \tilde{\pi}_i =
\begin{pmatrix} 1 & 0 \\ 0 & 0 \end{pmatrix}, \quad i=1,2.$$
Substitute $\tilde{\pi}_i$ into $J_2$ to get
\begin{eqnarray*} J_2 &=& J_0
(e^{-i\theta}\tilde{\pi} +e^{i\theta} \tilde{\pi}^\perp) (
-e^{i\theta}\tilde{\pi}_2
-e^{-i\theta} \tilde{\pi}_2^\perp) \\
&=& J_0 \left( \begin{pmatrix} e^{-i\theta} & 0 \\ 0 & e^{i\theta}
\end{pmatrix} +\xi \begin{pmatrix} 0 & -2i\sin\theta\
\bar{f}_1 \\ -2i\sin\theta\ f_1 & 0
\end{pmatrix} +O(\xi^2) \right) \\ & & \times
\left(
\begin{pmatrix} -e^{i\theta}
& 0 \\ 0 & -e^{-i\theta} \end{pmatrix} +\xi \begin{pmatrix} 0
& -2i \sin\o\ \bar{\beta}_1\bar{f}_1 \\
-2i \sin\o\ \beta_1 f_1 & 0
\end{pmatrix} \right. \\ & & \quad \left. +\xi \begin{pmatrix} 0
& -2i\sin\o\ \bar{\beta}_2\bar{f}_2 \\
-2i\sin\o\ \beta_2 f_2 & 0
\end{pmatrix} +O(\xi^2) \right)  \\ &=& J_0 \left( -I +
\xi \begin{pmatrix} 0 & c_1\bar{f}_1 \\ -\bar{c}_1  f_1 & 0
\end{pmatrix} +\xi \begin{pmatrix} 0 & c_2\bar{f}_2 \\ -\bar{c}_2
f_2 & 0 \end{pmatrix} +O(\xi^2) \right),
\end{eqnarray*} where $c_1, c_2\in {\mathbb C}$ are constants.
It follows that $\lim_{t\to -\infty} J_2 =-J_0$.
Note that $$\xi \begin{pmatrix} 0 & c_i\bar{f}_i \\
-\bar{c}_i f_i & 0 \end{pmatrix}, \quad i=1,2$$ is equal to the
unstable mode $\eta^+_{j_i,l_i}(c_i)$ at $-J_0$ given in
Proposition \ref{ak}, where $(j_1,l_1)=(j,l)$ and
$(j_2,l_2)=(j,2m-l)$. In other words, we have shown
\begin{equation} \lim_{t\to -\infty} \left( J_2 +J_0 + J_0\sum_{i=1}^2
\eta^+_{j_i,l_i}(c_i)\right) =0.\label{neginf} \end{equation} To
analyze the asymptotic behavior of $J_2$ as $t\to +\infty$, we set
$$\rho =e^{-\sqrt{m^2 - (j-m)^2
-(l-m)^2}\ t},\ h_1 =e^{i(jx+ly+mt)},\ h_2
=e^{i(jx+(2m-l)y+mt)}.$$ Then $\lim_{t\to +\infty} \rho=0$. A
similar computation implies that

\begin{itemize}
 \item[(1)] $\tilde{q}_1$ is parallel to $\begin{pmatrix}
 \rho h_1 \\ 1  \end{pmatrix}$.
 \item[(2)] $$\tilde{\pi}_1 =\begin{pmatrix} 0 & 0 \\ 0 & 1
 \end{pmatrix} +\rho \begin{pmatrix} 0 & h_1 \\ \bar{h}_1 & 0
 \end{pmatrix} +O(\rho^2),$$ and $$\tilde{\pi}_2  = \begin{pmatrix}
 0 & 0 \\ 0 & 1 \end{pmatrix} +\rho \begin{pmatrix} 0 &
 \gamma_1 h_1 \\ \bar{\gamma}_1 \bar{h}_1 & 0 \end{pmatrix}
 +\rho \begin{pmatrix} 0 & \gamma_2 h_2 \\
 \bar{\gamma}_2 \bar{h}_2 & 0 \end{pmatrix} +O(\rho^2)$$ for
 some constants $\gamma_1$, $\gamma_2$.
 \item[(3)] $$J_2 =J_0\left( -I +\rho \begin{pmatrix} 0 &
 d_1 h_1 \\ -\bar{d}_1 \bar{h}_1 & 0 \end{pmatrix}
 +\rho \begin{pmatrix} 0 & d_2 h_2 \\ -\bar{d}_2
 \bar{h}_2 & 0 \end{pmatrix} +O(\rho^2) \right)$$
 for some constants $d_1$, $d_2$. It follows that
 $\lim_{t\to +\infty} J_2 =-J_0$.
 \item[(4)] \begin{equation}\lim_{t\to +\infty} \left( J_2+J_0 +J_0
 \sum_{i=1}^{2} \eta_{j_i,l_i}^- (d_i) \right)
 =0.\label{posinf} \end{equation}
\end{itemize}

Formulas (\ref{neginf}) and (\ref{posinf}) imply that $J_2: T^2
\times {\mathbb R} \to SU(2)$ is a homoclinic Ward map.

Applying B\"{a}cklund transformations even times, with pairs of
poles and Hermitian projections chosen as above, we obtain more
homoclinic Ward maps. The case for $m<0$ is similar.  We
summarize the above discussion to give:

\begin{thm}
Let $m$ be a nonzero integer, $a=\text{diag} (im,-im)$, and $J_0
=e^{-(x+y)a}$. Choose $(j_{2k-1}, l_{2k-1}) \in {\mathbb Z}^2$
such that $$(j_{2k-1}-m)^2 +(l_{2k-1}-m)^2  <m^2, \quad
m<l_{2k-1}$$ and $(j_{2k-1}, l_{2k-1}) \not= (j_{2h-1}, l_{2h-1})$
for $1\le k<h \le N$. Let $(j_{2k}, l_{2k}) =(j_{2k-1},
2m-l_{2k-1})$,  $z_s =\sqrt{\frac{2m-j_s}{j_s}}\ e^{i\theta_s}$
with $\cos \theta_s =\sgn(m)\frac{(l_s-m)}{\sqrt{j_s (2m-j_s)}}$,
and $\sin \theta_s =\sgn(m)\sqrt{1-\cos^2\theta}$, $s=1,\dots,2N$.
Let $\pi_s$ be the Hermitian projection onto ${\mathbb C}
\begin{pmatrix} 1
\\ f_s(w_s)\end{pmatrix},$ where $f_s(w_s)
=e^{i(2m-j_s -2m z_s)w_s}$ and $w_s =x+z_su +z_s^{-1} v$.  Let
\begin{equation} J_{2N} (x,y,t) =h_{z_{2N},\pi_{2N}}*(\cdots
*(h_{z_1,\pi_1}*J_0 (x,y,t))\cdots )\end{equation} be the Ward map
obtained by applying $2N$ B\"{a}cklund transformations to $J_0$.
Then $J_{2N}: T^2 \times {\mathbb R} \to SU(2)$ is a homoclinic
Ward map. Moreover, $J_{2N} \to (-1)^N J_0$ as $t \to \pm\infty$.
\end{thm}

\begin{proof}
We have shown the $N=1$ case.  For general $N$, we use induction
and the calculation is similar. \end{proof}

The above construction can be generalized to $SU(n)$ model easily:

\begin{cor} Let $m$, $p$ be integers, $m\not= 0$, $1\le p\le n-1$,
$$a=\bpm i(n-p)mI_p & 0\\ 0 & -ipmI_{n-p} \epm \in su(n),$$
$\psi =e^{((1-\lambda)x +(1+\lambda-\lambda^2)u -v)a}$ the
extended solution, and $J_0=e^{-(x+y)a}$ the associated Ward map.
Choose $(j_{2k-1}, l_{2k-1}) \in {\mathbb Z}^2$ such that
$$\left( j_{2k-1}-\frac{nm}{2}\right)^2 +\left(
l_{2k-1}-\frac{nm}{2}\right)^2  < \left( \frac{nm}{2}\right)^2,
\quad \frac{nm}{2}<l_{2k-1}$$ and $(j_{2k-1}, l_{2k-1}) \not=
(j_{2h-1}, l_{2h-1})$ for $1\le k<h \le N$. Let $(j_{2k}, l_{2k})
=(j_{2k-1}, nm-l_{2k-1})$,  $z_s =\sqrt{\frac{nm-j_s}{j_s}}
e^{i\theta_s}$ with $\cos \theta_s =\sgn(m)
\frac{(l_s-nm/2)}{\sqrt{j_s (nm-j_s)}}$, and $\sin \theta_s
=\sgn(m) \sqrt{1-\cos^2\theta}$, $s=1,\dots,2N$. Let
$\pi_s(x,y,t)$ be the Hermitian projection of ${\mathbb C}^n$ onto
$${\mathbb C} (1,\cdots,1,f_s(w_s),\cdots f_s(w_s))^T,$$
where $1$ is repeated $p$-times, $f_s(w_s) =e^{i(nm-j_s -nm
z_s)w_s}$ is repeated $(n-p)$-times, $w_s =x+z_su +z_s^{-1} v$.
Let
\begin{equation} J_{2N} (x,y,t) =h_{z_{2N},\pi_{2N}}*(\cdots
*(h_{z_1,\pi_1}*J_0 (x,y,t))\cdots )\end{equation} be the Ward map
obtained by applying $2N$ B\"{a}cklund transformations to $J_0$.
Then $J_{2N}: T^2 \times {\mathbb R} \to SU(n)$ is a homoclinic
Ward map. Moreover, $J_{2N} \to (-1)^N J_0$ as $t \to \pm\infty$.
\end{cor}

The method discussed above can also produce Ward maps, which are
homoclinic to (time) periodic orbits. There are only some minor
changes in the construction, so we just list the main steps for
the $SU(2)$ model.

\begin{itemize}
  \item Let $m>0$ be an integer, and
  $b={\rm diag\/} (im,-im)$. Then
  $$\psi =e^{(x+(\l+2)u)b}$$ is an extended solution, and the
  associated Ward map is
  $$J_0=\psi^{-1}\vert_{\l=0}=e^{-(x+2u)b} =e^{-(x+y+t)b},$$
  which is triply periodic in the variables $x,y,t$.
  \item Set $\eta=J^{-1} \delta J$. Then the linearization of the
  Ward equation at $J_0$ is
  $$\eta_{tt}- \eta_{xx} -\eta_{yy} +[b, \eta_x +2\eta_y -2\eta_t]
  =0.$$ The unstable subspace of the linearization of the Ward
  equation at $J_0$ is $\bigoplus W_{jl}^+$, where $(j,l) \in
  {\mathbb Z}^2$, $(j-m)^2 +(l-2m)^2  <m^2$, and $W_{jl}^+ $
  is spanned by
  $$e^{\sqrt{m^2-(j-m)^2 -(l-2m)^2}\ t} \begin{pmatrix}
  0 & ce^{i(jx+ly+2mt)} \\ -\bar{c} e^{-i(jx+ly+2mt)}
  & 0 \end{pmatrix}$$ with constant $c\in \C$.
  The stable subspace at $J_0$ is $\bigoplus W_{jl}^-$,
  where $(j,l)$ satisfies the same condition, and $ W_{jl}^-$
  is spanned by
  $$e^{-\sqrt{m^2-(j-m)^2 -(l-2m)^2}\ t} \begin{pmatrix} 0 &
  ce^{i(jx+ly+2mt)} \\ -\bar{c} e^{-i(jx+ly+2mt)} & 0 \end{pmatrix}$$
  with constant $c\in \C$.
  \item Choose $(j,l)\in {\Z}^2$ with $(j-m)^2 +(l-2m)^2  <m^2$.
  Apply B\"{a}cklund transformation $h_{z_1,\pi_1}*\psi$, where
  $z_1 =re^{i\theta}$ with $r=\sqrt{\frac{2m-j}{j}}$, $\cos\theta
  =\frac{l-2m}{\sqrt{j(2m-j)}}$, $\sin\theta
  >0$, $\pi_1(x,y,t)$ is the Hermitian
  projection of ${\C}^2$ onto $\C \bpm 1 \\ e^{2im\a_1 w_1}
  \epm$, $\a_1 =\frac{2m-j}{2m}$, and $w_1=x+z_1 u+z_1^{-1}v$. Then
  $$\psi_1 = h_{z_1,\pi_1}*\psi =h_{z_1,\ti{\pi}_1}\psi,$$
  where $\ti \pi_1$ is the projection onto
  $$\psi(z_1)\Im\pi_1=\C \bpm 1 \\ e^{\sqrt{m^2-(j-m)^2 -(l-2m)^2}\,
  t} e^{-i(jx+ly+2mt)} \epm.$$
  \item Choose $(j_2,l_2)=(j,4m-l)\in\Z^2$. Apply B\"{a}cklund
  transformation again to get
  $$\psi_2 = h_{z_2,\pi_2}*\psi_1 =h_{z_2,\ti{\pi}_2}\psi_1.$$
  Here $z_2=-\bar{z}_1$, and $\pi_2(x,y,t)$ is the Hermitian projection
  onto $\C q$, where $q= \bpm 1 \\ e^{2im\a_2 w_2} \epm$, $\a_2 =\a_1$,
  and $w_2=x-\bar{z}_1 u -\bar{z}_1^{-1}v$. Then $\ti
  \pi_2(x,y,t)$ is the projection onto
  $$\C h_{z_1,\ti{\pi}_1}
  (-\bar{z}_1) \bpm 1 \\ e^{\sqrt{m^2-(j-m)^2 -(l-2m)^2}\, t}
  e^{-i(jx+(4m-l)y+2mt)} \epm.$$
  \item
  $$J_2=\psi_2^{-1} \vert_{\l=0} =J_0 (e^{-i\theta}
  \ti{\pi}_1 + e^{i\theta} \ti{\pi}_1^\perp) (-e^{i\theta}
  \ti{\pi}_2  -e^{-i\theta} \ti{\pi}_2^\perp)$$ is a Ward map
  from $T^2 \times\R$ to $SU(2)$. Analyzing the asymptotic
  behavior of $J_2$ as $t \to\pm\infty$, we see that
  $J_2$ is transversal and homoclinic to the periodic orbit
  $-J_0$. Applying B\"{a}cklund transformations even times with pairs
  of poles and Hermitian projections chosen similarly, we obtain
  more Ward maps which are homoclinic to $\pm J_0$.
\end{itemize}
The construction of homoclinic orbits to (time) periodic solutions
for the $SU(n)$ model is similar. Thus we have

\begin{thm} Let $m$, $p$ be integers, $m\not= 0$, $1\le p\le n-1$,
$$b=\bpm i(n-p)mI_p & 0\\ 0 & -ipmI_{n-p} \epm \in su(n),$$
$\psi =e^{(x+(\l+2)u)b}$ the extended solution, and $J_0
=e^{-(x+y+t)b}$ the associated Ward map. Choose $(j_{2k-1},
l_{2k-1}) \in {\mathbb Z}^2$ such that
$$\left( j_{2k-1}-\frac{nm}{2}\right)^2 +\left(
l_{2k-1}-nm \right)^2  < \left( \frac{nm}{2}\right)^2, \quad
nm<l_{2k-1}$$ and $(j_{2k-1}, l_{2k-1}) \not= (j_{2h-1},
l_{2h-1})$ for $1\le k<h \le N$. Let $(j_{2k}, l_{2k}) =(j_{2k-1},
2nm-l_{2k-1})$,  $z_s =\sqrt{\frac{nm-j_s}{j_s}} e^{i\theta_s}$
with $\cos \theta_s =\sgn(m)\frac{(l_s-nm)}{\sqrt{j_s (nm-j_s)}}$,
and $\sin \theta_s =\sgn(m)\sqrt{1-\cos^2\theta}$, $s=1,\dots,2N$.
Let $\pi_s(x,y,t)$ be the Hermitian projection of ${\mathbb C}^n$
onto
$${\mathbb C} (1,\cdots,1,f_s(w_s),\cdots f_s(w_s))^T,$$
where $1$ is repeated $p$-times, $f_s(w_s) =e^{i(nm-j_s)w_s}$ is
repeated $(n-p)$-times, and $w_s =x+z_su +z_s^{-1} v$. Let
\begin{equation*} J_{2N} (x,y,t) =h_{z_{2N},\pi_{2N}}*(\cdots
*(h_{z_1,\pi_1}*J_0 (x,y,t))\cdots )\end{equation*} be the Ward map
obtained by applying $2N$ B\"{a}cklund transformations to $J_0$.
Then $J_{2N}: T^2 \times {\mathbb R} \to SU(n)$ is a homoclinic
Ward map. Moreover, $J_{2N} \to (-1)^N J_0$ as $t \to \pm\infty$.
\end{thm}

\end{document}